\documentclass[12pt]{amsart}
\usepackage{amsthm,amsbsy,amsfonts,amssymb,amsmath,amscd}
\usepackage{latexsym}
\usepackage{fontenc}
\usepackage[mathscr]{euscript}

\input pictex.tex

\font\thinlinefont=cmr5

\theoremstyle{plain}
\newtheorem{thm}{Theorem}

\newtheorem{cor}[thm]{Corollary}

\theoremstyle{definition}
\newtheorem{defn}[thm]{Definition}

\newtheorem{rmk}[thm]{Remark}

\numberwithin{thm}{section}
\numberwithin{equation}{section}

\newcommand{\ga}[2]{\begin{gather}\label{#1}#2 \end{gather}}

\newcommand{\Spec}{{\rm Spec \,}}

\newcommand{\Gal}{{\rm Gal}}

\newcommand{\righthookdown}{\raise10pt%
  \hbox{$\scriptscriptstyle\cap$}%
  \kern-8.45pt\big\downarrow}

\newcommand{\dbQ}{{\mathbb{Q}}}
\newcommand{\dbZ}{{\mathbb{Z}}}
\newcommand{\Kbar}{\bar{K}}
\newcommand{\kbar}{\bar{k}}
\newcommand{\Vbar}{\bar{V}}
\newcommand{\Ybar}{\bar{Y}}
\newcommand{\GL}{\mathrm{GL}}
\newcommand{\scrC}{{\mathscr{C}}} 
\newcommand{\scrL}{{\mathscr{L}}} 
 
\newcommand{\Q}{{\mathbb Q}}

\newcommand{\Z}{{\mathbb Z}}

\begin{document}
\title[integrality]{Appendix
}
\author{Pierre Deligne}
\address{The Institute of Advanced Study, School of  Mathematics, 
NJ 08540 Princeton, USA}
\email{deligne@math.ias.edu}
\author{H\'el\`ene Esnault}
\address{
Universit\"at Duisburg-Essen, FB6, Mathematik, 45117 Essen, Germany}
\email{esnault@uni-essen.de}
\thanks{Partially supported by the 
DFG-Schwerpunkt ``Komplexe Mannigfaltigkeiten'' and by the DFG Leibniz Preis.}

\date{September 16, 2004}

\maketitle
 We generalize in this appendix Theorem 1.5 to nontrivial coefficients on varieties $V$ which are neither smooth nor projective. 
 We thank Alexander Beilinson,  Luc Illusie and Takeshi Saito
for very helpful discussions.\\ \ \\

The notations are as in the article. Thus $K$ is a local field with finite residue field $k$, $R\subset K$ is the ring of integers, $\Phi$ is a lifting of the geometric Frobenius in the Galois group of $K$. We consider
  $\ell$-adic sheaves on schemes of finite type defined over $K$ in the sense of \cite{
DeWeII}, (1.1). One generalizes the definition \cite{DeInt}, D\'efinition 5.1 
of $T$-integral $\ell$-adic  sheaves on schemes of finite type
 defined over finite fields to $\ell$-adic sheaves  on schemes of finite
type  defined over local fields with finite residue field. Recall to this aim that if $\scrC$ is a $\ell$-adic sheaf on a $K$-scheme $V$ 
and $v$ is a closed point of $V$, then 
the stalk $\scrC_{\bar{v}}$ of $\scrC$ at $\bar{v}$ is a ${\rm Gal}(\bar{K}/K_v)$-module, where $K_v\supset K$ is the residue field of $v$,
 with residue field $\kappa(v)\supset k$.   On $\scrC_{\bar{v}}$ the inertia $I_v={\rm Ker} ({\rm Gal}(\bar{K}/K_v)\to {\rm Gal}(\kappa(v)/k))$  acts quasi-unipotently (\cite{Gr2}).
 Consequently the eigenvalues
of a lifting $\Phi_v\in {\rm Gal}(\bar{K}/K_v) $ of the geometric Frobenius $F_v \in {\rm Gal}
(\overline{\kappa(v)}/\kappa(v))$ are, 
up to multiplication by roots of unity, well defined (\cite{DeWeII}, Lemma (1.7.4)). 
 Let $T\subset \Z$ be a set  of  prime numbers.
\begin{defn} \label{defn:int} The $\ell$-adic sheaf $\scrC$ is $T$-integral if the eigenvalues of $\Phi_v$ acting on $\scrC_{\bar{v}}$
are integral over $\Z[\frac{1}{t}, t\in T]$ for all closed points $v\in V$.
\end{defn}

\begin{thm} \label{thm:int}
Let $V$ be a scheme of finite type  defined over  $K$, and let $\scrC$ be a $T$-integral $\ell$-adic sheaf on $V$. Then if $f: V\to W$ is a morphism to another $K$-scheme of finite type  $W$ defined over $K$, the $\ell$-adic  sheaves 
$R^if_!\scrC$ are $T$-integral as well.
More precisely, if $w\in W$ is a closed point, then both 
$F_w$ and $|\kappa(w)|^{n-i}F_w$ acting on $(R^if_!\scrC)_{\bar{w}}$ 
  are integral over $\Z[\frac{1}{t}, t\in T]$, with $n={\rm dim}(f^{-1}(w))$.
\end{thm}
\begin{proof}
Let $w$ be a closed point of $W$.
By base change for $Rf_!$ and by
$\{w\}\hookrightarrow W$, one is reduced to the case
where $W$ is the spectrum of a finite extension $K'$
of $K$.
If $\Vbar:=V\otimes_{K'}\Kbar'$, for $\Kbar'$ an
algebraic closure of $K'$, one has to check the
integrality statements for the eigenvalues of a
lifting of Frobenius on $H_c^i(\Vbar,\scrC)$.

Let us perform the same reductions as in \cite{DeInt} p.~24.
Note that in loc. cit. $U$ can be shrunk so as to be
affine, with the sheaf smooth on it($=$ a local
system).
This reduces us to the cases where $V$ is of
dimension zero or is an affine irreducible curve
smooth over $W=\Spec(K')$ and $\scrC$ is a smooth 
sheaf.

The integrality statement to be proven is insensitive
to a finite extension of scalars $K''/K'$.
The $0$-dimensional case reduces in this way to the
trivial case where $V$ is a sum of copies of
$\Spec(K')$.
In the affine curve case, $H_c^0$ vanishes, while, as
in \cite{DeInt} Lemma 5.2.1, there is a $0$-dimensional
$Z\subset V$ such that the natural
($\Phi$-equivariant) map from $H^0(Z,\scrC)(-1)$ to
$H_c^2(\bar{V},\scrC)$ is surjective, leaving us
only $H_c^1$ to consider.

Let $\scrC_{\dbZ_\ell}$ be a smooth $\dbZ_\ell$-sheaf
from which $\scrC$ is deduced by $\otimes\dbQ_\ell$,
and let $\scrC_\ell$ be the reduction modulo $\ell$
of $\scrC_{\dbZ_\ell}$.
For some $r$, it is locally (for the \'etale topology)
isomorphic to $(\dbZ/\ell)^r$.
Let $\pi\colon\, V'\to V$ be the \'etale covering of
$V$ representing the isomorphisms of $\scrC_\ell$
with $(\dbZ/\ell)^r$.
It is a $\GL(r,\dbZ/\ell)$-torsor over $V$.
As $H_c^*(\Vbar,\scrC)$ injects into
$H_c^*(\overline{V'},\pi^*\scrC)$, renaming irreducible
components of $V'$ as $V$ and $K'$ as $K$, we may and
shall assume that $W=\Spec(K)$ and that $\scrC_\ell$
is a constant sheaf.

Let $V_1$ be the projective and smooth completion of
$V$, and $Z:=V_1\setminus V$.
Extending scalars, we may and shall assume that $Z$
consists of rational points and that $V_1$, marked
with those points, has semi-stable reduction.
It hence is the general fiber of $X$ regular and
proper over $\Spec(R)$, smooth over
$\Spec(R)$ except for quadratic non-degenerate
singular points, with $Z$ defined by disjoint
sections $z_\alpha$ through the smooth locus.

Let $Y$ be the special fiber of $X$, and $\Ybar$ be
$Y\times_k\kbar$, for $\kbar$ the residue field of
the algebraic closure $\Kbar$ of $K$.

$$
\begin{matrix}
V &\null\kern5.0 true cm
\hidewidth:\text{\small{complement of 
     disjoint sections $z_\alpha$}}\hidewidth &&&\\
\raise3pt\hbox{$\scriptstyle j$}
   \righthookdown &&&&\\
V_1 &\hookrightarrow &X &\hookleftarrow &Y\\
\big\downarrow 
&&\big\downarrow &&\big\downarrow\\
{\text{\small Spec}}(K) &\hookrightarrow
&{\text{\small Spec}}(R) &\hookleftarrow 
  &{\text{\small Spec}}(k) \end{matrix}
$$
The cohomology with compact support
$H_c^1(\Vbar,\scrC)$ is  $H^1(\overline{V_1},j_!\scrC)$, and
vanishing cycles theory relates this $H^1$ to the
cohomology groups on $\Ybar$ of the nearby cycle
sheaves $\psi^i(j_!\scrC)$, which are $\ell$-adic
sheaves on $\Ybar$, with an action of $\Gal(\Kbar/K)$
compatible with the action of $\Gal(\Kbar/K)$
(through $\Gal(\kbar/k)$) on $\Ybar$.
The choice of a lifting of Frobenius, i.e. of a
lifting of $\Gal(\kbar/k)$ 
in $\Gal(\Kbar/K)$, makes them come
from $\ell$-adic sheaves on $Y$, to which the
integrality results of \cite{DeInt} apply.
Using the exact sequence
$$
0\to H^1(\Ybar,\psi^0(j_!\scrC))\to
H^1(\overline{V_1},j_!\scrC)\to H^0(\Ybar,\psi^1(j_!\scrC))
$$
and \cite{DeInt} Th\'eor\`eme 5.2.2, we are reduced to check integrality of
the sheaves $\psi^i(j_!\scrC)$ ($i=0,1$).
It even suffices to check it at any $k$-point $y$ of
$Y$, provided we do so after any unramified finite
extension of $K$.

Let $X_{(y)}$ be the henselization of $X$ at $y$, and
$Y_{(y)}$, $V_{1(y)}$ and $V_{(y)}$ be the inverse
image of $Y$, $V_1$ or $V$ in $X_{(y)}$.
There are three cases:
$$
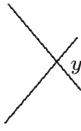
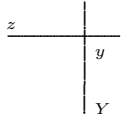
\vbox{\beginpicture
\setcoordinatesystem units <.70cm,.70cm>
\unitlength=1.00000cm
\setshadesymbol ({\thinlinefont .})
\setlinear
%
%
\linethickness= 0.500pt
\setplotsymbol ({\thinlinefont .})
{\plot  2.667 22.822  4.032 21.203 /
}%
%
%
\linethickness= 0.500pt
\setplotsymbol ({\thinlinefont .})
{\plot  2.599 20.574  3.965 22.193 /
}%
%
%
\put{$\scriptstyle y$}%
 [lB] at  3.842 21.584
%
%
\linethickness= 0.500pt
\setplotsymbol ({\thinlinefont .})
{\plot 18.603 20.828 18.610 22.947 /
}%
%
%
\put{$\scriptscriptstyle\bullet$}%
 [lB] at 18.510 21.827
%
%
\put{$\scriptstyle y$}%
 [lB] at 18.986 21.795
%
%
\linethickness= 0.500pt
\setplotsymbol ({\thinlinefont .})
{\plot 11.779 20.750 11.786 22.868 /
}%
%
%
\linethickness= 0.500pt
\setplotsymbol ({\thinlinefont .})
{\plot 12.446 22.221 10.331 22.227 /
}%
%
%
\put{$\scriptscriptstyle z$}%
 [lB] at 10.289 22.352
%
%
\put{$\scriptscriptstyle Y$}%
 [lB] at 11.972 20.733
%
%
\put{$\scriptscriptstyle y$}%
 [lB] at 11.972 21.844
\linethickness=0pt
\putrectangle corners at  2.574 22.972 and 18.986 20.549
\endpicture}
$$

\noindent
\parbox{6.5 true in}{(1) $y$ singular on $Y$
\kern2.0cm (2) $y$ on a $z_\alpha$
\kern2.0cm (3) general case.}

\bigskip\noindent
The restriction of $\psi^i$ to $y$ depends only on
the restriction of $\scrC$ to $V_{(y)}$, and short 
exact sequences of sheaves give rise to long exact
sequences of $\psi$.

Because $\scrC_\ell$ is a constant sheaf, $\scrC$ is
tamely ramified along $Y$ and the $z_\alpha$. 
More precisely, it is
given by a representation of the pro-$\ell$
fundamental group of $V_{(y)}$.
It is easier to describe the group deduced from the
profinite fundamental group by pro-$\ell$ completing
only the kernel of its map to
$\widehat{\dbZ}=\Gal(\kbar/k)$.
By Abhyankhar's lemma, this group is an extension of
$\widehat{\dbZ}$, generated by Frobenius, by
$\dbZ_\ell(1)^2$ in case (1) or (2) or $\dbZ_\ell(1)$
in case (3).
The representation is given by $r\times r$ matrices
congruent to $1\mod\ell$.
For $\ell\not=2$, such a matrix, if quasi-unipotent,
is unipotent.
Indeed, it is the exponential of its logarithm and
the eigenvalues of its logarithm are all zero.
For $\ell=2$, the same holds if the congruence is
$\mod4$, hence if  $\scrC_{\dbZ_2}\mod4$ is constant,
a case to which one reduces by the same argument we
used $\mod2$.
By Grothendieck's argument \cite{SeTa} p.515, the action of
$\dbZ_\ell(1)$ or $\dbZ_\ell(1)^2$ is
quasi-unipotent, hence unipotent, and we can filter
$\scrC$ on $V_{(y)}$ by smooth sheaves such that
the successive quotients $Q$ extend to 
smooth sheaves on $X_{(y)}$.
If $Q$ extends to a smooth sheaf $\scrL$ on
$X_{(y)}$, the corresponding $\psi$ are known by
Picard-Lefschetz theory: $\psi^0$ is $\scrL$ restricted to
$Y_{(y)}$ in cases (1) and (3), and $\scrL$ outside
of $y$ extended by zero in case (2); $\psi^1$ is 
non-zero only in case (1), where it $\scrL(-1)$ on
$\{y\}$ extended by zero.
By d\'evissage, this gives the required integrality.
\end{proof}
\begin{cor}\label{cor:int} Let $V$ be smooth scheme of finite type defined 
over $K$.  Then the eigenvalues of $\Phi$ on $H^i(\bar{V}, \Q_\ell)$ are integral over $\Z$. 
\end{cor}
\begin{proof}
If $K$ has characteristic zero, there is a good compactification
$j:V\hookrightarrow W$, with $W$ smooth proper over $K$ and $D=W\setminus V=
\cup D_i$ a
strict normal crossing divisor. Then the long exact sequence
\ga{a.5}{
\ldots \to H^i_D(\bar{W}, \Q_\ell)\to H^i(\bar{W}, \Q_\ell) \to 
H^i(\bar{V}, \Q_\ell)\to \ldots}
and Theorem \ref{thm:int} applied to the cohomology of $W$ reduces to showing integrality for $H^i_D(\bar{W}, \Q_\ell)$. As in (3.3) of the article, the Mayer-Vietoris spectral sequence
\ga{a.6}{E_1^{-a+1,b}=
\oplus_{|I|=a} H^{b}_{D_I}(\bar{W}, \Q_\ell)
\Rightarrow H^{1-a+b}_D(\bar{W}, \Q_\ell),}
with $D_I=\cap_{i\in I} D_i$, 
reduces to the case where $D$ is smooth projective of codimension $\ge 1$.
Then  
 purity together with
 Theorem \ref{thm:int}
allow to conclude. If $K$ has equal positive characteristic, we apply de Jong's theorem \cite{deJ}, Theorem 6.5 to find $\pi: V'\to V$ generically finite and
$j:V'\hookrightarrow W$ a good compactification. As $\pi^*: 
H^i(\bar{V}, \Q_\ell)\hookrightarrow H^i(\overline{V'}, \Q_\ell)$ is injective, we conclude as above. 
\end{proof}
Corollary \ref{cor:int} gives some flexibility as we do not assume 
that $V$ is projective. In particular, 
one can apply the same argument as in the proof of Theorem 2.1 of the article in order to show an improved version of Theorem 1.5, (ii) there:
\begin{cor} \label{cor:int2}
Let $V$ be a smooth scheme of finite type over $K$, and $A\subset V$ be a 
codimension $\kappa$ subscheme. Then the eigenvalues of $\Phi$ on $H^i_A(\bar{V}, \Q_\ell)$ are divisible by $|k|^\kappa$ as algebraic integers.
\end{cor}
\begin{proof}One has a stratification $\ldots \subset A_i
\subset A_{i-i} \subset \ldots A_0=A$ by closed subschemes defined over $K$ with $A_{i-1}\setminus A_i$ smooth.  
The $\Phi$-equivariant long exact sequence
\ga{a.7}{
\ldots \to H^m_{A_i}(\bar{V}, \Q_\ell)\to H^m_{A_{i-1}}(\bar{V}, \Q_\ell)
\to  H^m_{(A_{i-1}\setminus A_i)}(\overline{V\setminus A_i}, \Q_\ell)
\to \ldots
}
together with purity and  Corollary \ref{cor:int} allow to conclude
 by induction on the codimension. 
\end{proof}
\begin{rmk}
One has to pay attention that even if Theorem \ref{thm:int}  generalizes
 Theorem 1.5 i) of the article to $V$ not necessarily smooth, there is no such generalization of  Theorem 1.5 ii) to the non-smooth case, even on a finite field. Indeed, let $V$ be a  rational curve 
with one node. Then $H^1(\bar{V}, \Q_\ell)=\Q_\ell(0)$ as we see from the normalization sequence, yet $H^1_{{\rm node}}(\bar{V}, \Q_\ell)=
H^1(\bar{V}, \Q_\ell)$ as the localization map $H^1(\bar{V}, \Q_\ell)\to 
H^1(\overline{V \setminus {\rm node}}, \Q_\ell) $ 
factorizes through $H^1(\overline{{\rm normalization}}, \Q_\ell)=0$.  
So we can't improve the integrality statement to a divisibility statement in 
general. In order to force divisibility, one needs the divisor supporting the cohomology to be in good position with respect to the singularities. 

\end{rmk}

\bibliographystyle{plain}

\begin{thebibliography}{99}
\bibitem{deJ} de Jong, A. J.: Smoothness, semi-stability and alterations, Publ. Math. IHES {\bf 83} (1996), 51-93.
\bibitem{DeInt} Deligne, P.: Th\'eor\`eme d'int\'egralit\'e, Appendix  to
 Katz, N.: Le niveau de la cohomologie des intersections compl\`etes, 
Expos\'e XXI
in SGA 7, Lect. Notes Math. vol. {\bf 340}, 363-400,
Berlin Heidelberg New York Springer 1973.
\bibitem{DeWeII} Deligne, P.: La conjecture de Weil, II. Publ. Math. IHES 
{\bf 52} (1981), 137-252.
\bibitem{Gr2} Grothendieck, A.: Groupes de monodromie en g\'eom\'etrie 
alg\'ebrique, SGA 7 I,
Lecture Notes in Mathematics {\bf 288}, Springer Verlag. 
\bibitem{SGA1} Raynaud, Mme: Propret\'e cohomologique des 
faisceaux d'ensembles et des faisceaux de groupes non commutatifs, 
in Rev\^etements \'Etales et Groupes Fondamentaux, SGA 1, expos\'e XIII, Lecture Notes in Mathematics {\bf 224}, Springer Verlag. 
\bibitem{SeTa} Serre, J.-P.; Tate, J.: Good reduction of abelian varieties, 
Annals of Mathematics {\bf 88} (1968), 492-517.
\end{thebibliography}

\renewcommand\refname{References}

\end{document}